\input amstex
\documentstyle{amsppt}

\magnification=1200
\parindent 20pt

\vsize=7.5in
\NoBlackBoxes
\TagsOnRight

\define \a{\alpha}
\define \be{\beta}
\define \Dl{\Delta}
\define \dl{\delta}

\define \s{\sigma}

\define \ve{\varepsilon}

\define \fc{\frac}
\define \iy{\infty}

\define \ode{\overset{\de}\to{=}}

\define \ov{\overline}

\define \edm{\enddemo}
\define \ep{\endproclaim}

\define \bk{\bigskip}

\define \1{^{-1}}
\define \2{^{-2}}

\define \CA{\Cal A}

\define \cD{\Cal D}

\define \CL{\Cal L}

\define  \CK{\Cal K }
\define  \CT{\Cal T }

\define \loc{\operatorname{loc}}
\define \de{\operatorname{def}}

\define \xir{x\in R}
 \baselineskip 15pt
 \topmatter

\heading{\bf{Some properties of the Sturm-Liouville operator in
$L_p(R)$}}\endheading
\vskip .05in

 \centerline{N.A. Chernyavskaya}
\centerline{Department of Mathematics and Computer Science}
\centerline{Ben-Gurion University of the Negev}
\centerline{P.O.B. 653, Beer-Sheva, 84105, Israel}
\vskip .05in
\centerline{L.A. Shuster}
\centerline{Department of Mathematics and Computer Science}
\centerline{Bar-Ilan University}
\centerline{Ramat-Gan, 52900, Israel}

\abstract  We consider the boundary problem
$$- y''(x)+q(x)y(x)=f(x),\quad\ x\in R\tag1$$
$$\lim_{|x|\to\iy}y^{(i)}(x)=0,\quad i=0,1\tag2$$
where $f(x)\in L_p(R),\ p\in[1,\iy],$ \ $1\le q(x)\in L_1^{\loc}(R).$
For this boundary problem we obtain: 1) necessary and sufficient conditions
for unique solvability and a priori
properties of the solution; 2) a criterion for the resolvent to be compact in
$ L_p(R),$\ $p\in [1,\iy]$ and some a priori properties of the
spectrum.\endabstract\endtopmatter

\document
\baselineskip 20pt
\subheading{\S1. Introduction}

In this paper we consider the boundary problem (1.1) - (1.2):
$$- y''(x)+q(x)y(x)=f(x),\quad\ x\in R\tag1.1$$
$$\lim_{|x|\to\iy}y^{(i)}(x)=0,\quad i=0,1\tag1.2$$
where $f(x)\in L_p(R),\ p\in[1,\iy]$\ $(\|\cdot\|_\iy$ is defined as
essup), and $$1\le q(x)\in
L_1^{\loc}(R).\tag1.3$$

Our aim is to give a detailed description of the Sturm-Liouville operator
$\CL_p$ in $L_p(R),$\ $p\in[1,\iy]$ (see
below).
Usually, the operator $\CL_p$ is defined as the closure in $L_p(R)$ of a
differential  expression $\ell:$
$$\ell y=-y''(x)+q(x)y(x)\tag1.4$$
with $y(x)$ belonging to $C_0^\iy(R)$ \cite{10, Ch.VII, \S5]}, \cite{1,
Ch.II, \S1}.
Such an approach allows one to develop a sufficiently complete theory of
the operator $\CL_p$, but some problems
are left in the background or not considered at all.
This "neglect" is quite natural since the above approach is built upon
general functional methods rather than on
concrete properties of the differential operator $\ell.$
Below we add some new facts to already known ones using a more detailed
study of the operator $\ell$ and its
inversion.

\proclaim{Theorem  1.1}

\rom{A)} \ For $p\in[1,\iy)$ for any $f(x)\in L_p(R)$ , the boundary
problem $(1.1)$-$(1.2)$ has a unique solution
$y(x)\in\cD_p,$ where
$$\cD_p=\{z(x): z(x)\in C^{(1)}(R)\cap L_p(R),\ z'(x)\in  AC^{\loc}(R),
(\ell z)(x)\in L_p(R), $$
$$\lim_{|x|\to\iy}z^{(i)}(x)=0,\quad i=0,1\}$$

\rom{B)}\ Let $p=\iy.$
For any $f(x)\in L_\iy(R)$, the boundary problem $(1.1)$-$(1.2)$ has a
solution (moreover, a unique solution) if
and only if $(1.5)$ holds:
$$\lim_{|x|\to\iy}\int_{x-a}^{x+a}q(t)dt=\iy\quad \text{for any}\
a\in(0,\iy).\tag1.5$$
In the latter case, the solution $y(x)$ to $(1.1)$-$(1.2)$ belongs to
$\cD_\iy^0,$ where
$$\align\cD_\iy^0&=\{z(x): z(x)\in C^{(1)}(R),\ z'(x)\in A C^{\loc}(R),\
(\ell z)(x)\in L_\iy(R),\\ & \lim_{|x|\to \iy}
z^{(i)}(x)=0,\quad i=0,1\}.\endalign$$

\rom{C)}\ Let $p=\iy,$ and suppose that $(1.5)$ fails.
For any $f(x)\in L_\iy(R)$, equation $(1.1)$ has a unique solution $y(x)\in
\cD_\iy$ where
$$\cD_\iy=\{z(x):z(x)\in C^{(1)}(R),\ z'(x)\in A C^{\loc}(R),\ (\ell
z)(x)\in L_\iy(R)\}.$$\ep

In the sequel, by the Sturm-Liouville operator $\CL_p$, we mean the
differential expression (1.4) defined on
$\cD_p.$ Thus $\CL_p:\cD_p\to L_p(R),\ p\in [1,\iy],$ and there exists
$\CL_p\1,$\ $\cD(\CL_p\1)=L_p(R),$\
$p\in[1,\iy]$ (see Theorem 1.1).
The following fact is well known for $p=2$ (see \cite{8}).
\proclaim{Theorem 1.2} The operator $\CL_p\1:L_p(R)\to L_p(R)$ is bounded
for $p\in[1,\iy].$
Moreover, for any $p\in[1,\iy]$ it is compact if and only if $(1.5)$ holds.
In the latter case, the spectrum $\s(\CL_p)$ of the operator $\CL_p$ is
purely discrete for all $p\in[1,\iy]$ and
does not depend on $p$ for $p\in[1,\iy)$.
In particular, $\s(\CL_p)=\s(\CL_2)\subset [1,\iy)$ for $p\in[1,\iy).$\ep
\remark{{\bf Concluding remarks and acknowledgements}} Some of our results
have been proved in \cite{13} under the
additional assumption
$$1\le q(x)\in L_p^{\loc} (R),\quad p\in[1,\iy].\tag1.6$$
This paper was stimulated by crucial remarks of the referee concerning our
paper \cite{4}. In \cite{4} we used
results of \cite{13}, and the referee noticed that they can be made more
precise.
In particular, he suggested the definition of the classes $\cD_p,\cD_\iy^0$
as ``a natural $L_p$-framed for
(1.1)" (we cite the report).
Theorem 1.2 (under condition (1.6)) was a subject of the authors'
discussion with Prof. V.G. Maz'ya in 1989.
It was V.G. Maz'ya who noticed that Theorem 1.2 must hold under the weakest
requirement (1.5) replacing (1.6).
The authors thank the referee of \cite{4} and V.G. Maz'ya for their
attention, benevolent criticism, and
formulation of the problems.
\bk

\subheading{\S2. \ Preliminaries}

Throughout the sequel we denote by $c$ absolute positive constants whose
values are not essential for exposition
and may differ even within a single chain of calculations.
\proclaim{Lemma 2.1}
Consider an equation
$$z''(x)=q(x)z(x),\quad \xir.\tag2.1$$
There exists a fundamental system of solutions (FSS)
$\{u(x),v(x)\}$ of (2.1) such that
$$\aligned &u(x)>0,\ v(x)>0,\ u'(x)<0,\ v'(x)>0,\quad \xir\\
&v'(x)u(x)-u'(x)v(x)=1,\ u(x)=v(x)\int_x^\iy\fc{dt}{v^2(t)},\quad \xir\\
&\lim_{x\to\iy}u(x)=\lim_{x\to\iy}u'(x)=\lim_{x\to-\iy}
v(x)=\lim_{x\to-\iy}v'(x)=0
\\&
\lim_{x\to\iy}v(x)=\lim_{x\to\iy}v'(x)=\lim_{x\to-\iy}u(x)=\lim_{x\to-\iy}
|u'(x)|=\iy.\endaligned\tag2.2$$\ep
Under condition (1.3), equation (2.1) does not  oscillate at $\pm\iy,$ and
therefore $u(x),v(x)$ are determined up
to constant factors as principal solutions to (2.1) on $(0,\iy)$ and
$(-\iy,0)$, respectively \cite{7,p.355}.
We call a FSS of (2.1) with   properties (2.2) a principal FSS of (2.1)
(PFSS, \cite{3}).
The proof of Lemma 2.1 can be also found in \cite{3}.
\proclaim{Lemma 2.2} \cite{5} For $\xir,$ a PFSS of $(2.1)$ admits the
following representation:
$$v(x)=\sqrt{\rho(x)}\exp\left(\fc{1}{2}\int_{x_0}^x\fc{dt}{\rho(t)}\right),\
u(x)=\sqrt{\rho(x)}\exp\left(-\fc{1}{2}\int_{x_0}^x\fc{dt}{\rho(t)}\right),\ 
\xir\tag2.3$$
where $\rho(x)\ode u(x)v(x)$, $\xir$  and $x_0$ is the unique root of the
equation $u(x)=v(x).$

Moroever $$|\rho'(x)|<1,\quad\xir.\tag2.4$$\ep
Formulae of type (2.3) were also mentioned in \cite{14, pp.419-420}.
  Lemma 2.2 can be found in the above presentation in \cite{3}.

For a fixed $\xir,$ consider the following equations in $d\ge0:$
$$1=\int_0^{\sqrt 2d} \int_{x-t}^xq(\xi)d\xi dt,\quad 1=\int_0^{\sqrt
2d}\int_x^{x+t}q(\xi)d\xi dt,\quad
2=d\int_{x-d}^{x+d}q(t)dt.\tag2.5$$
For any $\xir,$ each of the equations (2.5) has a unique finite positive
solution \cite{3}.
Denote the solutions by $d_1(x),\ d_2(x),\ d(x)$ respectively.
\proclaim{Lemma 2.3} The function $d(x)$ is continuous for $\xir$ and
$$0<d(x)\le 1,\quad \xir.\tag2.6$$
Moreover, for $\ve\in[0,1],$  \ $\xir$ the following inequalities hold:
$$(1-\ve)d(x)\le d(t)\le (1+\ve)d(x),\quad t\in[x-\ve d(x),x+\ve
d(x)].\tag2.7$$\ep
Lemma 2.3 and the definition of $d(x)$ are due to M.O. Otelbaev (see \cite{9}).
The functions $d_1(x),\ d_2(x)$ were introduced in \cite{2}.
Lemma 2.3 was proved in the above presentation in \cite{4}.
\proclaim{Theorem 2.1 \cite{3}} For $\xir$ one has the following estimates
(see $(2.3)$):
$$\fc{1}{\sqrt 2}\le\fc{v'(x)}{v(x)}d_1(x)\le \sqrt 2,\quad \fc{1}{\sqrt
2}\le\fc{|u'(x)|}{u(x)}d_2(x)\le\sqrt
2,\tag2.8$$
$$\fc{1}{\sqrt 2}\fc{d_1(x)d_2(x)}{d_1(x)+d_2(x)}\le \rho(x)\le\sqrt
2\fc{d_1(x)d_2(x)}{d_1(x)+d_2(x)},\tag2.9$$
$$4\1 d(x)\le\rho(x)\le 3\cdot 2\1 d(x).\tag2.10$$
\ep
\proclaim{Lemma 2.4} \cite{13} Condition $(1.5)$ holds if and only if
$\lim\limits_{|x|\to\iy}d(x)=0.$
\ep
\proclaim{Theorem 2.2} \cite{6} Let $p\in[1,\iy).$
A set $\CK\subset L_p(R)$ is precompact if and only if the following
conditions hold:

$1)$\ $\sup\limits_{f\in\CK}\|f\|_p<\iy;$

$2)$\
$\lim\limits_{\dl\to0}\sup\limits_{f\in\CK}\sup\limits_{|t|\le\dl}\|f(\cdot
+t)-f(\cdot)\|_p=0$
\quad $3)$\ $\lim\limits_{N\to\iy}\sup\limits_{f\in\CK}\int_{|x|\ge
N}|f(x)|^pdx=0.$\ep
\proclaim{Theorem 2.3} \cite{11} Let $X,$\ $Y$ be normed spacaes, $T: X\to
Y$ a linear bounded operator, $T^*$ the
adjoint operator.
Then $T$ is compact if and only if $T^*$ is compact.\ep

\proclaim{Theorem 2.4} \cite{12} Let $(X,\Sigma,\mu)$ be a space with a
$\Sigma$-infinite measure, $L_0$ the space
of piecewise constant $\mu$-integrable functions on $X,$ \ $M$ the space of
all measure functions on $X,$\ $T:
L_0\to M$ a linear operator having continuous extensions $T_p: L_p(\mu)\to
L_p(\mu),$\ $T_1: L_q(\mu)\to L_q(\mu)$
where $1\le p\le q\le \iy.$
Let the spectra $\s(T_p)$ and $\s(T_q)$ be zero-dimensional (say, countable).
Then for all $r\in [p,q]$ one has $\s(T_r)=\s(T_q).$
(The operator $T$ has a continuous extension $T_r: L_r(\mu)\to L_r(\mu),\
r\in[p,q]$ by the Riesz-Torin Theorem
(see \cite{6, Ch.VI, \S10}.)\ep
\bk

\subheading{\S3.\ Unconditional solvability of   boundary problem in
$L_p(R),\ p\in [1,\iy)$}

In this section, we prove part A) of Theorem 1.1
We need Lemmas 3.1 -- 3.4.

\proclaim{Lemma 3.1} For a PFSS of $(2.1)$ one has the following inequalities:
$$v'(x)\ge v(x),\quad |u'(x)|\ge u(x),\quad \xir.\tag3.1$$
$$v(x)\ge\exp(x-t)v(t)\quad\text{for}\quad x\ge t;\
u(x)\ge\exp(t-x)u(t)\quad\text{for}\quad x\le t.\tag3.2$$
$$\rho(x)=u(x)v(x)\le 1,\quad\xir.\tag3.3$$
\ep
\demo{Proof} Let $\hat z(t)=\exp(t).$
Then for $t\in R$ we obtain
$$[v'(t)\hat z(t)-\hat z'(t)v(t)]'=v''(t)\hat z(t)-\hat
z''(t)v(t)=(q(t)-1)v(t)\hat z(t).\tag3.4$$
Taking into account (1.3) and (2.2), we deduce inequality (3.1) from (3.4)
for $v(x):$
$$\hat z(x)[v'(x)-v(x)]=v'(x)\hat z(x)-\hat
z'(x)v(x)=\int_{-\iy}^x(q(t)-1)v(t)\hat z(t)dt\ge 0,\quad \xir.$$
Together with (2.2), this implies (3.2) for $v(\cdot):$
$$ln\fc{v(x)}{v(t)}=\int_t^x\fc{v'(\xi)}{v(\xi)}d\xi\ge x-t,\quad x\ge t.$$
The following chain of calculations is based on (2.2) and (3.1) and leads
to (3.3):
$$\align
\rho(x)&=v^2(x)\int_x^\iy\fc{d\xi}{v^2(\xi)}=-v^2(x)\int_x^\iy\fc{1}{v'(\xi)}
d\left(\fc{1}{v(\xi)}\right)
\\
&=-v^2(x)\left[\fc{1}{v'(\xi)v(\xi)}\left|_x^\iy\right.+\int_x^\iy
 \fc{1}{v(\xi)}
 \fc{v''(\xi)d\xi}{v'(\xi)^2}\right]\\
&=\fc{v(x)}{v'(x)}-v^2(x)\int_x^\iy\fc{q(\xi)d\xi}{v'(\xi)^2}\le
\fc{v(x)}{v'(x)}\le 1.\endalign$$
The proof of (3.1) -- (3.2) for $u(\cdot)$ is similar to the above proof
for $v(\cdot).$\hfill $\qed$\edm
\smallskip
Let us introduce the Green function $G(x,t)$ and the Green operator $(Gf)(x):$
$$G(x,t)=\cases u(x)v(t),\quad & x\ge t\\
u(t)v(x),\quad & x\le t\endcases\tag3.5$$
$$(Gf)(x)=\int_{-\iy}^\iy G(x,t)f(t)dt,\quad \xir,\quad f(\cdot)\in L_p(R)$$
\proclaim{Lemma 3.2} For $x,t\in R$ one has the inequalities
$$0<G(x,t)\le\exp(-|t-x|),\quad 0<G(x,t)\le\fc{3}{4}d(x)\exp(-|t-x|).\tag3.6$$
$$\left|\fc{\partial}{\partial x}G(x,t)\right|\le\exp(-|t-x|).\tag3.7$$
\ep
\demo{Proof} Let $x\ge t.$
Then (3.5), (3.3), and (3.2) imply (3.6):
$$G(x,t)=u(x)v(t)=\rho(x)\fc{v(t)}{v(x)}\le \fc{v(t)}{v(x)}\le \exp(-|t-x|).$$
The case $x\le t$  can be considered in a similar way.
If in both cases ($x\ge t$ and $x\le t$) we use (2.10) instead of (3.3), we
obtain the second inequality of (3.6).
Let $x>t.$
Then (3.5), (2.2), and (3.2) imply (3.7):
$$\left|\fc{\partial}{\partial x} G(x,t)\right|=|u'(x)|v(t)
=|u'(x)|v(x)\fc{v(t)}{v(x)}\le
\fc{v(t)}{v(x)} \le\exp(-|t-x|).$$
The case $x<t$ can be treated in a similar way; the case $x=t$ follows from
(2.4).
\edm
\proclaim{Lemma 3.3}
Let $r(x)=\fc{1}{d(x)},\ \xir.$
The functions $(Gf)(x),\ \fc{d}{dx}(Gf)(x),$ where $f(x)\in L_p(R),\
p\in[1,\iy]$ are absolutely continuous and
satisfy the inequalities:
$$\|r(x)(Gf)(x)\|_p\le c\|f\|_p,\quad \|r(x)(Gf)(x)\|_{C(R)}\le
c\|f\|_p.\tag3.8$$
$$\left\|\fc{d}{dx}(Gf)(x)\right\|_p\le c\|f\|_p,\quad
\left\|\fc{d}{dx}(Gf)(x)\right\|_{C(R)}\le
c\|f\|_p.\tag3.9$$
\ep
\demo{Proof} For $f(x)\in L_p(R),\ p\in[1,\iy]$ consider the integrals
$$\CT_1(x)=\int_{-\iy}^x v(t)f(t)dt,\quad \CT_2(x)=\int_x^\iy
u(t)f(t)dt,\quad\xir.\tag3.10$$
These integrals converge absolutely since (3.2) and H\"older's inequality
imply:
$$|\CT_1(x)|\le\int_{-\iy}^x v(t)|f(t)|dt\le
v(x)\int_{-\iy}^x\exp(-|t-x|)|f(t)|dt\le cv(x)\|f\|_p.\tag3.11$$
$$|\CT_2(x)|\le\int_x^\iy u(t)|f(t)|dt\le
u(x)\int_x^\iy\exp(-|t-x|)|f(t)|dt\le cu(x)\|f\|_p.\tag3.12$$
Since $(Gf)(x)=u(x)\CT_1(x)+v(x)\CT_2(x),$ we conclude that both $(Gf)(x)$ and
$(Gf)'(x)=u'(x)\CT_1(x)+v'(x)\CT_2(x)$ are absolutely continuous.
The following estimates can be derived from (3.6) and H\"older's inequality:
$$\aligned
r(x)&|(Gf)(x)|\le r(x)\int_{-\iy}^\iy G(x,t)|f(t)|dt\le c\int_{-\iy}^\iy
\exp (-|t-x|)|f(t)|dt
\\
&\le
c\left(\int_{-\iy}^\iy\exp(-|t-x|)dt\right)^{1/p'}
\left(\int_{-\iy}^\iy\exp(-|t-x|)|f(t)|^pdt\right)^{1/p}\\
&\le c\left(\int_{-\iy}^\iy\exp(-|t-x|)|f(t)|^pdt\right)^{1/p}\le
c\|f\|_p\endaligned\tag3.13$$
In particular, (3.13) and (2.6) imply the following inequalities:
$$\|(Gf)(x)\|_{C(R)}\le \|r(x)(Gf)(x)\|_{C(R)}\le c\|f\|_p,\quad
p\in[1,\iy].\tag3.14$$
{}From (3.13), (2.6) and Fubini's theorem we get
$$\aligned \|G(f)(x)\|_p^p&\le\|r(x)(Gf)(x)\|_p^p\le
c\int_{-\iy}^\iy\left[\int_{-\iy}^\iy|f(t)|^p\exp(-|t-x|)dt\right]dx\\
&=c\int_{-\iy}^\iy|f(t)|^p\left[\int_{-\iy}^\iy\exp(-|t-x|)dx\right]dt\le
c\|f\|_p.\endaligned\tag3.15$$
We thus obtained estimates (3.8).
Inequalities (3.9) can be found in a similar way by applying (3.7) instead
of (3.6).\hfill$\qed$\edm
\proclaim{Corollary 3.3.1}
Let $f(x)\in L_p(R),\ p\in[1,\iy].$
The function $y(x)=(Gf)(x),\ \xir$ is the unique solution to $(1.1)$ in the
class $C^{(1)}(R).$
In particular, part $c)$ of Theorem $1.1$ holds.\ep
\demo{Proof}
It is an immediate consequence of Lemma 3.3 and (2.2)\hfill$\qed$\edm
\proclaim{Lemma 3.4}
Let $f(x)\in L_p(R),\ p\in[1,\iy).$
Then
$$\lim\limits_{|x|\to\iy}(Gf)(x)=0,\quad
\lim_{|x|\to\iy}\fc{d}{dx}(Gf)(x)=0.\tag3.16$$\ep
\demo{Proof} Let $\CA\in (0,\iy).$
{}From (3.13) and (2.6) $\Rightarrow$ (3.17), \quad (3.17) $\Rightarrow$ 
(3.18):
$$\aligned  |  (Gf)(x)|&\le c\int_{-\iy}^\iy |f(t)|^p\exp(-|t-x|)dt\le
c\int_{|t-x|\le\CA}|f(t)|^p\exp(-|t-x|)dt\\
&+c\int_{|t-x|\ge \CA}|f(t)|^p\exp(-|t-x|)dt\le c\int_{|t-x|\le
\CA}|f(t)|^pdt+c\exp(-\CA)\|f\|_p^p.\endaligned\tag3.17$$
$$0\le
\varliminf_{|x|\to\iy}|(Gf)(x)|^p\le\varlimsup_{|x|\to\iy}|(Gf)(x)|^p\le
c\exp(-\CA)\|f\|_p^p.\tag3.18$$
In (3.18) we pass to limit as $\CA\to\iy$ and obtain the first equality of
(3.16).
The second equality of (3.16) is checked similarly using (3.7).\hfill
$\qed$\edm

Theorem 1.1 A) follows from Lemma 3.1, Corollary 3.3.1 and (3.16).\hfill 
$\qed$\bk
\subheading{\S4.\ A criterion for solvability of a singular boundary
problem in $L_\iy(R)$}

In this section we prove Part B) of Theorem 1.1.

\demo{Necessity} Let $y(x)=(Gf)(x),\ f(x)\in L_\iy(R).$
{}From Corollary 3.3.1 it follows that the general solution of (1.1) has the
form $z(x)=\a u(x)+\be v(x)+y(x),\ \xir,\
\a,\be-\text{const.}$ From (2.2) and (3.8) we conclude that if $z(x)$ is a
solutoin to (1.1) -- (1.2), then
$\a=\be=0$ and $z(x)\equiv y(x),\ \xir.$
Let now $f(x)\equiv 1,\ \xir.$
In this case, the solution to (1.1) -- (1.2) is of the form
$$  y(x) =u(x)\int_{-\iy}^xv(t)dt+v(x)\int_x^\iy u(t)dt \ode
u(x)I_1(x)+v(x)I_2(x). $$ From (2.2) we easily derive lower estimates for
$$\align  I_1(x)&=\int_{-\iy}^xv(t)dt=\int_{-\iy}^x
v(t)[v'(t)u(t)-u'(t)v(t)]dt\ge\int_{-\iy}^x v(t)v'(t)u(t)dt\\
&\ge u(x)\int_{-\iy}^xv(x)v'(t)dt=\fc{u(x)v^2(x)}{2},\quad\xir,\endalign$$
$$\align  I_2(x)&=\int_{x}^\iy u(t)dt=\int_{x}^\iy
u(t)[v'(t)u(t)-u'(t)v(t)]dt\ge -\int_{x}^\iy v(t)u'(t)u(t)dt\\
&\ge -v(x)\int_{x}^\iy u'(x)u(t)dt=\fc{v(x)u^2(x)}{2},\quad\xir,\endalign$$
Hence $y(x)\ge (u(x)v(x))^2=\rho^2(x).$
Then (1.2), (2.10) and Lemma 2.4 imply (1.5). \edm

\demo{Sufficiency} By Corollary 3.3.1, the function $y(x)=(Gf)(x),\ f(x)\in
L_\iy(R)$ satisfies (1.1) almost
everywhere, $y(x)\in C^{(1)}(R),\  y'(x)\in A C^{\loc}(R);$
and by (3.6)
$$|y(x)|\le\int_{-\iy}^\iy G(x,t)|f(t)|dt\le
cd(x)\int_{-\iy}^\iy\exp(-|t-x|)dt\|f\|_\iy\le cd(x)\|f\|_\iy,\ \xir.$$
Hence by Lemma 2.4, $y(x)\to0$ as $|x|\to\iy.$

\proclaim{Lemma 4.1} Condition $(1.5)$ is equivalent to either of the
equalities $(4.1)$
$$\lim_{|x|\to\iy}d_1(x)=0,\quad \lim_{|x|\to\iy}d_2(x)=0.\tag4.1$$\ep
\demo{Proof}
Let us verify that (1.5) is equivalent to the first equality of (4.1).
(For the second equality of (4.1) this equivalence can be established
similarly.)
Suppose that (1.5) holds.
Note that from (1.3) it follows that
$d_1(x)\le 1$ for $\xir:$
$$1=\int_0^{\sqrt 2 d_1(x)}\int_{x-t}^xq(\xi)d\xi dt\ge\int_0^{\sqrt 2
d_1(x)}\int_{x-t}^x d\xi
dt=d_1^2(x)\Rightarrow 0<d_1(x)\le 1.\tag4.2$$
Let $a=\fc{2\sqrt 2}{3},$ \ $b=\fc{\sqrt 2}{3},$\ $\overline x=x-bd_1(x).$
Then (see \cite{3}):
$$\align
2&=2\int_0^{\sqrt 2 d_1(x)}\int_{x-t}^x q(\xi)d\xi dt\ge
2\int_{ad_1(x)}^{\sqrt 2 d_1(x)}\int_{x-t}^x q(\xi)d\xi
dt\\
&\ge 2bd_1(x)\int_{x-2bd_1(x)}^x q(\xi)d\xi\ge (bd_1(x))\int_{\overline
x-(bd_1(x))}^{\overline x+(bd_1(x))}
q(\xi)d\xi.\endalign$$
Hence $d(\ov x)\ge bd_1(x),\ \xir,$ i.e. $$3.2^{-1/2}d(x-\sqrt 2\cdot 3\1
d_1(x))\ge d_1(x)>0,\quad\xir.\tag4.3$$
Then Lemma 2.4, (4.2) and (4.3) imply (4.1).
Conversely, suppose that (4.1) holds. By the definition of $d_1(x),$
 $$2=2\int_0^{\sqrt 2 d_1(x)}\int_{x-t}^x q(\xi)d\xi dt\le 2\sqrt 2
d_1(x)\int_{x-\sqrt 2 d_1(x)}^x q(\xi)d\xi\le
(2\sqrt 2 d_1(x))\int_{x-(2\sqrt 2 d_1(x))}^{x+(2\sqrt 2 d_1(x))}.$$
Hence $2\sqrt 2 d_1(x)\ge d(x)>0$ for $\xir.$
Therefore $d(x)\to0$ as $|x|\to\iy,$ and it remains to use Lemma 4.2.\hfill
$\qed$\edm
\proclaim{Lemma 4.2} Suppose that $(1.5)$ holds. Then
$$\lim_{|x|\to\iy}\int_{-\iy}^\iy\left|\fc{\partial}{\partial x}
G(x,t)\right|dt=0.\tag4.4$$\ep
\demo{Proof}
{}From $(3.5),$ $(2.2),$ $(2.8),$ and $(2.9)$ we obtain for $\xir:$
$$\aligned \int_{-\iy}^\iy\left|\fc{\partial}{\partial x}
G(x,t)\right|dt&=|u'(x)\int_{-\iy}^xv(t)dt+v'(x)\int_x^\iy u(t)dt \\
&=\fc{|u'(x)|}{u(x)}\ \fc{\rho(x)}{v(x)}\int_{-\iy}^x
v(t)dt+\fc{v'(x)}{v(x)}\ \fc{\rho(x)}{u(x)}\int_x^\iy u(t)dt\\
&\le 2\left\{\fc{1}{v(x)}\int_{-\iy}^x v(t)dt+\fc{1}{u(x)}\int_x^\iy
u(t)dt\right\}.\endaligned\tag4.5$$
The equalities
$$\lim_{|x|\to\iy}\fc{1}{v(x)}\int_{-\iy}^x v(t)dt=0,\quad
\lim_{|x|\to\iy}\fc{1}{u(x)}\int_x^\iy u(t)dt=0.\tag4.6$$
are verified in the same way, and therefore we only check the first one.
{}From (2.2) and (3.1) we obtain:
$$v(x)\to\iy\quad\text{as}\quad x\to\iy\Rightarrow \int_{-\iy}^x v(t)dt\to
\iy \quad\text{as}\quad x\to\iy.\tag4.7$$
$$v'(x)\ge v(x)\Rightarrow v(x)\ge\int_{-\iy}^x v(t)dt,\ \xir\Rightarrow
\int_{-\iy}^x
v(t)dt\to0\ \text{as}\  x\to\iy.\tag4.8$$
{}From Lemmas 2.4 and 4.1 it follows that (4.1) holds, and therefore by (2.8),
$\lim\limits_{|x|\to\iy}\fc{v(x)}{v'(x)}=0.$
By L'\^opital's rule, taking into account (4.7) and (4.8), we get
$$\lim_{|x|\to\iy}\fc{\int_{-\iy}^x v(t)dt}{v(x)}
=\lim_{|x|\to\iy}\fc{v(x)}{v'(x)}=0.\tag4.9$$
Equality (4.4) now follows from (4.9) and (4.5).\hfill $\qed$\edm

To end the proof of the theorem, it remains to show that $y'(x)\to0$ as
$|x|\to\iy$ and to prove that (1.1) --
(1.2) has a unique solution. From (4.4) it follows that
$$|y'(x)|\le\int_{-\iy}^\iy\left|\fc{\partial}{\partial x} G(x,t)\right|\
|f(t)|dt\le\int_{-\iy}^\iy\left|\fc{\partial}{\partial
x}G(x,t)\right|dt\cdot\|f\|_\iy\to0\ \text{as}\ |x|\to\iy$$
and therefore $y(x)\in\cD_\iy^0.$
The uniqueness in the class $\cD_\iy^0$ follows from (2.2).\hfill$\qed$\edm
\bk
\subheading{\S5. \ Properties of the Green operator in $L_p(R)$}

This section is devoted to Theorem 1.2.
{}From Theorem 1.1 and the definition of the operator $\CL_p$ (see \S1) it
follows that $\CL_p\1=G$ (see (3.5))
By (2.6) and (3.8), $\|G\|_{p\to p}\le c<\iy,$\ $p\in[1,\iy].$
To verify that the condition of the theorem on compactness of $G$ is
necessary, we need the following Lemmas 5.1 and
5.2.
\proclaim{Lemma 5.1}
For $\xir,$\ $t\in\left[x-\fc{d(x)}{2},x+\fc{d(x)}{2}\right]$ one has the
inequalities
$$c\1v(x)\le v(t)\le cv(t),\quad c\1u(x)\le u(t)\le cu(x).\tag5.1$$\ep
\demo{Proof}
{}From (2.3), (2.7) (for $\ve=\fc{1}{2}),$ and (2.10), we get
$$\align\fc{v(t)}{v(x)}&\le\sqrt{\fc{\rho(t)}{\rho(x)}}\exp\left(\fc{1}{2}
\left|\int_x^t\fc{d\xi}
{\rho(\xi)}\right|\right)=\sqrt
{\fc{\rho(t)}{d(t)}
\cdot\fc{d(t)}{d(x)}\cdot\fc{d(x)}{\rho(x)}}\exp\left(\fc{1}{2}\left|\int_x^t
\fc
{d(\xi)}{\rho(\xi)}\cdot
\fc{d(x)}{d(\xi)}\cdot\fc{d\xi}{d(\xi)}\right|\right)\\
&\le c.\endalign$$
$$\align\fc{v(t)}{v(x)}&\ge\sqrt{\fc{\rho(t)}{\rho(x)}}\exp\left(-\fc{1}{2}
\left|\int_x^t\fc{d\xi}{\rho(\xi)}\right|\right)=\sqrt
{\fc{\rho(t)}{d(t)}\cdot\fc{d(t)}{d(x)}\cdot\fc{d(x)}{\rho(x)}}
\exp\left(-\fc{1}{2}\left|\int_x^t\fc{d(\xi)}{\rho(\xi)}\cdot
\fc{d(x)}{d(\xi)}\cdot\fc{d\xi}{d(\xi)}\right|\right)\\
&\ge c\1.\endalign$$
Inequalities (5.1) for $u(\cdot)$ are checked similarly.\hfill $\qed$\edm
\proclaim{Lemma 5.2}
Let $f_x(\xi)$ be the characteristic function of the segment
$\Dl(x)=\left[x-\fc{d(x)}{2},x+\fc{d(x)}{2}\right],$\
$\xir.$
Then for $t\in\Dl(x)$ one has the inequalities
$$(Gf_x)(t)\ge c\1d^2(x),\quad\xir.\tag5.2$$
\ep
\demo{Proof}
For $t\in\Dl(x),$ by (5.1) and (2.10), we get
$$\align (Gf_x)(t)&=u(t)\int_{x-\fc{1}{2}d(x)}^t
v(\xi)d\xi+v(t)\int_t^{x+\fc{1}{2}d(x)}u(\xi)d\xi\\
&=\left[\fc{u(t)}{u(x)}\int_{x-\fc{1}{2}d(x)}^t\fc{v(\xi)}{v(x)}d\xi\right]
\rho(x)+\left[\fc{v(t)}{v(x)}\int_t^{x+\fc{d(x)}
{2}}\fc{u(\xi)}{u(x)}d\xi\right]\rho(x)\\
&\ge c\1 \rho(x)d(x)\ge c\1d^2(x).\hfill  \qquad
\qquad\qquad\qquad\qquad\qquad\qquad\qquad\qquad\qed
\endalign$$ \edm

Let $p\in[1,\iy)$ and $f_x(\xi),$\ $\xir$ be the functions from Lemma 5.2.
By (2.6), one has $\|f_x\|_p\le 1,$\ $\xir.$
If $G: L_p(R)\to L_p(R)$ is compact, then by Theorem 2.2 we get
$$\align
0&= \lim_{N\to\iy}\sup_{\|f\|_{p\le
1}}\left[\int_{-\iy}^{-N}|(Gf)(t)|^pdt+\int_N^\iy|(Gf)(t)^pdt\right]\\
&\ge \lim_{N\to\iy}\sup_{|x|\ge
N}\left[\int_{-\iy}^{-N}|(Gf_x)(t)|_p^pt+\int_N^\iy|(Gf_x)(t)|^pdt\right]\\
&\ge \lim_{N\to\iy}\sup_{x\le
-N}\left[\int_{x-\fc{d(x)}{2}}^x|(Gf_x)(t)|^pdt\right]+\lim_{N\to\iy}
\sup_{x\ge N}\left[
\int_x^{x+\fc{d(x)}{2}}|Gf_x)(t)|^pdt\right]\\
&\ge \lim_{N\to\iy}\sup_{|x|\ge N}c\1d(x)^{2p+1}.\endalign$$
Hence $\lim\limits_{|x|\to\iy} d(x)=0,$ and (1.5) holds by Lemma 2.4.
To prove that the condition of the theorem on compactness of $G$ is
sufficient, we need Theorem 5.1.
\proclaim{Theorem 5.1}
Let $p\in [1,\iy).$
The operator $G:L_p(R)\to L_p(R)$ is compact if
$$\lim_{|x|\to\iy}\int_{-\iy}^\iy G(x,t)dt=0.\tag5.3$$
\ep
We divide the proof of the theorem into several separate assertions.
\proclaim{Lemma 5.3}
Suppose that $(5.3)$ holds.
Then
$$\lim_{N\to\iy}\sup_{\xir}\left(\int_{-\iy}^{-N}G(x,t)dt\right)=0,\quad
\lim_{N\to\iy}\sup_{\xir}\left(\int_N^\iy
G(x,t)dt\right)=0.\tag5.4$$\ep
\demo{Proof}
Equalities (5.4) are checked in the same way.
Let us prove, say, the second one.
For given $x$ and $N,$ in the cases $x\le N$ and $x\ge N,$ respectively, we
establish inequalities proving (5.4):
$$\align H(x,N)&\ode \int_N^\iy G(x,t)dt=v(x)\int_N^\iy u(t)dt\le
v(N)\int_N^\iy u(t)dt \le\int_{-\iy}^\iy
G(N,t)dt\\
&\le \sup_{|x|\ge N} \int_{-\iy}^\iy G(x,t)dt\endalign$$
$$H(x,N)=u(x)\int_N^x v(t)dt +v(x)\int_x^\iy u(t)dt\le\int_{-\iy}^\iy
G(x,t)dt\le \sup_{|x|\ge N}\int_{-\iy}^\iy
G(x,t)dt.\hfill \qed    $$\edm
\proclaim{Lemma 5.4}
If (5.3) holds, then for $p\in[1,\iy)$ one has
$$\lim_{N\to\iy}\sup_{\|f\|_p\le 1}\int_{|x|\ge N}|(Gf)(x)|^pdx=0.\tag5.5$$\ep
\demo{Proof} By H\"older's inequality, Fubini's theorem and (3.6), we get
$$\align \sup_{\|f\|_p\le 1}&\int_{|x|\ge N}|(Gf)(x)|^pdx\\
&\le\sup_{\|f\|_p\le 1}\left\{\int_{|x|\ge
N}\left[\int_{-\iy}^{\iy}G(x,t)dt\right]^{p/p'}\left[\int_{-\iy}^\iy
G(x,t)|f(t)|^pdt\right]dx\right\}\\
&\le c\sup_{\|f\|_p\le 1}\int_{|x|\ge N}\left[\int_{-\iy}^\iy
G(x,t)|f(t)|^pdt\right]dx\\&\le
 c\sup_{\|f\|_p\le 1}
\int_{-\iy}^\iy |f(t)|^p\left[\int_{|x|\ge N} G(x,t)dx\right]dt
\le c\ \sup_{t\in R}\int_{|x|\ge N} G(x,t)dx.\endalign$$
The lemma follows from (5.3), (5.4) taking into account that $G(x,t)$ is
symmetric.  \edm
\proclaim{Lemma 5.5}
Let $p\in[1,\iy)$ and suppose that $(5.3)$ holds.
Denote $y(x)=(Gf)(x),$\ $f(x)\in L_p(R).$
Then
$$\lim_{\eta\to0}\sup_{\|f\|_p\le1}\|y(\cdot+\eta)-y(\cdot)\|_p=0.\tag5.6$$\ep
\demo{Proof} By (3.9), we get for $\|f\|_p\le 1:$
$$|y(x+\eta)-y(x)|\le \left|\int_x^{x+\eta}|y'(\xi)|d\xi\right|\le|\eta| \
\left\|\fc{d}{dx}G\right\|_{p\to
C(R)} \|f\|_p\le c|\eta|.$$
Let $|\eta|\le 1.$
Then
$$\int_{|x|\ge N}|y(x+\eta)-y(x)|^pdx\le 2^{p+1}\int_{|x|\ge N-1}|y(x)|^pdx.$$
By Lemma 5.4, for a given $\ve>0$ there is $N(\ve)\gg 1$ such that
$$2^{p+1}\sup_{\|f\|_p\le 1}\left(\int_{|x|\ge
N(\ve)-1}|y(x)|^pdx\right)\le\ve.$$
Therefore, taking into acocunt the above proved inequalities, we get
$$\align&\lim_{\eta\to0}\sup_{\|f\|_p\le1}\int_{-\iy}^\iy|y(x+\eta)-y(x)|^p
dx=\lim_{\eta\to0}
\left\{\int_{-N(\ve)}^{N(\ve)}
|y(x+\eta)-y(x)|^pdx\right.\\
 &\ +\left.\int_{|x|\ge
 N(\ve)}|y(x+\eta)-y(x)|^pdx\right\}\le\lim_{\eta\to0}\left\{cN(\ve)|\eta|^p
+2^{p+1}\int_{|x|\ge
 N(\ve)-1}|y(x)|^pdx\right\}\\
 &\ \le\ve.\endalign$$
Since $\ve$ is an arbitrary positive number, (5.6) is proved.\hfill $\qed$\edm
Since the operator $G: L_p(R)\to L_p(R)$ is bounded (see above), Theorem
5.1 follows from (5.6), (5.5) and Theorem
2.2.\hfill $\qed$
\medskip
To deduce the sufficiency of the condition of Theorem 1.2 on the
compactness of $G: L_p(R)\to L_p(R)$ for
$p\in[1,\iy),$ note that (3.6) implies
$$0\le\int_{-\iy}^\iy G(x,t)dt\le cd(x)\int_{-\iy}^\iy \exp(-|t-x|)dt\le
cd(x).\tag5.7$$
If (1.5) holds, then Lemma 2.4 and (5.7) imply (5.3), and by Theorem 5.1,
the operator $G: L_p(R)\to L_p(R)$ is
compact for $p\in[1,\iy).$
Since $G(x,t)=G(t,x),$\ $x,t\in R,$ we conclude that
$(\CL_1\1)^*=\CL_\iy\1$, and the operator $\CL_\iy\1:
L_\iy(R)\to L_\iy(R)$ is compact if and only if $\CL_1\1: L_1(R)\to L_1(R)$
is compact (see Theorem 2.3), i.e., if
and only if (1.5) holds.
Thus the criterion for compactness of $\CL_p\1: L_p(R)\to L_p(R)$,\
$p\in[1,\iy]$ is proved.
Let us check the other assertions of Theorem 1.1.
Under condition (1.5), for any $p\in[1,\iy]$ the spectrum $\s(\CL_p\1)$ is
at most countable, and for $p\in[1,\iy)$
the operator $\CL_p\1=G$ satisfies the conditions of Theorem 2.4.
Hence in that case $\s(\CL_p)=\s(\CL_2),$\ $p\in[1,\iy).$
The operator $\CL_2\1$ is symmetric and bounded, thus it is self-adjoint.
Hence $\s(\CL_2)\in R.$
Moreover, by (1.3) the operator $\CL_2$ is semi-bounded from below because
for $y(x)\in\cD(\CL_2)=\cD_2$ one has
$$\aligned \langle \CL_2y,y\rangle&=\int_{-\iy}^\iy (-y''(x)+q(x)y(x))y(x)dx\\
&=y'(x)y(x)\bigm|_{-\iy}^\iy+\int_{-\iy}^\iy[y'(x)^2+q(x)y^2(x)]dx\\&=
\int_{-\iy}^\iy[y'{}^2(x)+q(x)y^2(x)]dx\ge \int_{-\iy}^\iy
|y(x)|^2dx= \|y\|_2^2\endaligned\tag5.8$$ From (5.8) it follows that
$\s(\CL_2)\in[1,\iy).$
Theorem 1.2 is proved.
\hfill$\qed$
\bk

\end